\documentclass[12pt,a4paper]{article}
\usepackage{latexsym,fullpage,amsfonts,amssymb,amsmath,amscd,graphics,epic,theorem}
\usepackage[all]{xy}
\usepackage{amsmath}
\usepackage{amssymb, amstext}
\baselineskip 18pt \textwidth 15cm  \theoremstyle{plain}
\newtheorem{lemma}{Lemma}[section]
\newtheorem{proposition}[lemma]{Proposition}
\newtheorem{remark}[lemma]{Remark}

\newtheorem{theorem}[lemma]{Theorem}
\newtheorem{definition}[lemma]{Definition}
\newtheorem{notation}[lemma]{Notation}

\newtheorem{corollary}[lemma]{Corollary}

{\theorembodyfont{\rmfamily}

\newcommand{\C}{{\mathbb C}}

\newcommand{\F}{{\mathbb F}}
\newcommand{\E}{{\mathbb E}}
\newcommand{\proofend}{\hfill$\Box$\smallskip}

\newcommand{\Span}{{\operatorname{Span}}}

\newcommand{\Ker}{{\operatorname{Ker}}}

\newcommand{\Sc}{{\cal S}}

\newcommand{\tG}{{\widetilde{G}}}
\newcommand{\tO}{{\widetilde{O}}}
\newcommand{\Fou}{{\mathcal{F}}}

\begin{document}

\title{A proof of the multiplicity one conjecture\\ for $\mathrm{GL}_n$ in $\mathrm{GL}_{n+1}$}

\author{Avraham Aizenbud and Dmitry Gourevitch}
\maketitle
\begin{abstract}
Let $\mathbb{F}$ be a non-archimedean local field of characteristic
zero. We consider distributions on $\mathrm{GL}(n+1,\mathbb{F})$
which are invariant under the adjoint action of
$\mathrm{GL}(n,\mathbb{F})$. We prove that any such distribution is
invariant with respect to transposition. This implies that the
restriction to $\mathrm{GL}(n)$ of any irreducible smooth
representation of $\mathrm{GL}(n+1)$ is multiplicity free.

Our paper is based on the recent work \cite{RS} of Steve Rallis and
G\'{e}rard Schiffmann where they made a remarkable progress on this
problem.

In \cite{RS}, they also show that our result implies multiplicity 
one theorem for restrictions from the orthogonal group $\mathrm{O}(V
\oplus \mathbb{F})$ to $\mathrm{O}(V)$.
\end{abstract}
\setcounter{tocdepth}{1}
\subsection*{Acknowledgements}
We would like to thank our teacher {\bf Joseph Bernstein} for
teaching us most of the mathematics we know.

We cordially thank {\bf Joseph Bernstein} and {\bf Eitan Sayag} for
guiding us through this project. Their contribution was very
important.

We would also like to thank {\bf Vladimir Berkovich}, {\bf Yuval
Flicker}, {\bf Erez Lapid}, {\bf Omer Offen} and {\bf Yiannis
Sakellaridis} for useful remarks.
 \tableofcontents
\section{Introduction}
Let $\F$ be a non-archimedean local field of characteristic zero.
Consider the standard imbedding $\mathrm{GL}(n, \F) \subset
\mathrm{GL}(n+1, \F)$. Let $\mathrm{GL}(n, \F)$ act on
$\mathrm{GL}(n+1, \F)$ by conjugation. The goal of this paper is to
prove the following theorem.
\begin{theorem} \label{Goal}
Every $\mathrm{GL}(n, \F)$-invariant distribution on
$\mathrm{GL}(n+1, \F)$ is invariant with respect to transposition.
\end{theorem}
This theorem is important in representation theory, since it implies
the following multiplicity one theorem (see e.g. \cite{RS}, section
2).
\begin{theorem} \label{Rep1}
Let $\pi$ be an irreducible smooth representation of
$\mathrm{GL}(n+1,\F)$ and $\rho$ be an irreducible smooth
representation of $\mathrm{GL}(n,\F)$. Then $$\dim
\mathrm{Hom}_{\mathrm{GL}(n,\F)} (\pi|_{\mathrm{GL}(n,\F)},\rho)
\leq 1.$$
\end{theorem}
In their recent paper \cite{RS}, in part II, Steve Rallis and
G\'{e}rard Schiffmann have shown that theorem \ref{Goal} also
implies similar theorems for distributions on unitary and orthogonal
groups, which in turn imply multiplicity one results for those
groups.

In \cite{RS}, Rallis and Schiffmann have also made a remarkable
progress in proving the above theorem \ref{Goal}, and our paper is
based on their results. For the benefit of the reader we present the
proofs of all the statements from \cite{RS} that we use.

Theorem \ref{Goal} also gives another proof of Bernstein's theorem
about $P$-invariant distributions on $\mathrm{GL}(n)$ (see
\cite{Ber}) which proves Kirillov's conjecture in the
non-archimedean case.
\subsection{Reformulation of the main theorem}
Let $G:=G_n :=\mathrm{GL}(n,\F)$. Consider the action of the
2-element group $S_2$ on $G$ given by the involution $g \mapsto {\,
^tg^{-1}}$. It defines a semidirect product $\tG := \tG_n := G_n
\rtimes S_2$. Let $V:=V_n=\F^n$ and $X:=X_n := \mathrm{sl}(V_n)
\times V_n \times V_n^*$ where $\mathrm{sl}(V)\subset
\mathrm{End}(V)$ is the space of operators with zero trace.

The group $\tG$ acts on $X$ by
$$g(A,v,\phi):= (gAg^{-1},gv, {g^{-1}}^* \phi) \text{ and }$$
$$T(A,v,\phi):=({}^tA,{}^t\phi,{}^tv)$$
where $g \in G$ and $T$ is the generator of $S_2$. Here, ${}^tA$
denotes the transposed  matrix in $sl_n$, ${}^t\phi \in V_n$ denotes
the column vector corresponding to the row vector $\phi \in V_n^*$,
and ${}^tv$ denotes the row vector corresponding to the column
vector $v \in V_n$. Also for any operator $g:V \to V$, we denote by
$g^*:V^* \to V^*$ the adjoint operator.

Note that $\tG$ acts separately on $\mathrm{sl}(V)$ and on $V \times
V^*$. Define a character $\chi$ of $\tG$ by $\chi(g,s):=
\mathrm{sign}(s)$.

Theorem \ref{Goal} can be deduced from the following theorem.
\begin{theorem} \label{main}
Any $(\tG,\chi)$-equivariant distribution on $X$ is zero.
\end{theorem}
The deduction was first proven in \cite{RS} (section 5). We prove it
in section \ref{Reform} in a slightly different way. In section
\ref{Reform} we also give a coordinate-free definition of the group
$\tG$ and its action on $X$.

\subsection{Sketch of our proof}

The theorem will be proved by induction on $n$. Let $S$ denote the
closure of the union of the supports of all $(\tG,
\chi)$-equivariant distributions on $X$. We would like to prove that
$S = \emptyset$.

Let $pr_1:X \to \mathrm{sl}(V)$ and $pr_2:X \to V \oplus V^*$ be the
natural projections. Rallis and Schiffman have shown that the
induction hypothesis
implies:\\
(i) $pr_2(S)$ is contained in $Y:=\{(v,\phi)\in V \oplus V^*| \,
\langle \phi, v\rangle =0\}$\\
(ii) $pr_1(S)$ is contained in the nilpotent cone $\mathcal{N}$.\\
Part (i) follows from the localization principle and Frobenious reciprocity.\\
Part (ii) is proven using Harish-Chandra descent method.\\ We will
present a proof of this statement in the appendix (section
\ref{SecRalSchiff}).

For any vector $v$ and covector $\phi$ let $v \otimes \phi$ denote
the operator of rank one defined by them. Let $\lambda \in \F$ be a
scalar. We use a family of automorphisms $\nu_{\lambda}$ of $X$
defined by
$$\nu_{\lambda}(A,v,\phi):=(A+\lambda v\otimes \phi-\lambda \frac{\langle
\phi,v \rangle}{n}\mathrm{Id},v,\phi).$$  Note that the
automorphisms $\nu_{\lambda}$ commute with the action of $\tG$ and
hence preserve $S$.

Let $\mathcal{N}_i\subset \mathrm{sl}(V)$ be the union of all
nilpotent orbits of dimensions $\leq i$. We prove by downward
induction that $S \subset \mathcal{N}_i\times Y$ for all $i$.
Suppose $S \subset \mathcal{N}_i\times Y$. We have to show $S
\subset \mathcal{N}_{i-1}\times Y$. Since $\nu_{\lambda}(S)=S$, $S$
is contained in the intersection $\bigcap
\nu_{\lambda}(\mathcal{N}_i\times Y)$.

We have to show that for any nilpotent orbit $O$ of dimension $i$,
the restriction of any $(\tG, \chi)$-equivariant distribution $\xi$
to $O \times Y$ is zero. As we have seen, the support of
$\xi|_{O\times Y}$ is contained in $(O \times Y) \cap (\bigcap
\nu_{\lambda}(\mathcal{N}_i\times Y))$, which we denote by $\tO$.
Using the fact that the Fourier transform of a
$(\tG,\chi)$-equivariant distribution is also
$(\tG,\chi)$-equivariant, the theorem boils down to the following
key lemma.
\begin{lemma} [Key lemma]
Let $O$ be a nilpotent orbit.  Let $\zeta \in \Sc^*(O\times V \times
V^*)^{\tG,\chi}$. Suppose that both $\zeta$ and $\Fou(\zeta)$ are
supported in $\widetilde{O}$. Then $\zeta=0$.
\end{lemma}

Using Frobenious reciprocity, the key lemma reduces to a statement
about distributions on $V \oplus V^*$.

Namely, fix $A\in O$. Let $R_A$ denote the fiber over $A$ of the
projection $\widetilde{O} \to O$. Then $\zeta$ corresponds to a
distribution $\eta$ on $V \oplus V^*$
with the following properties:\\
(i)$\, $   $\eta$ is supported in $R_A$\\
(ii)$\,$  $\Fou(\eta)$ is supported in $R_A$\\
(iii) $\eta$ is $\chi$-equivariant with respect to the stabilizer of
$A$ in $\tG$.

We have to show $\eta=0$. We will prove that $R_A$ is contained in
$$Q_A :=\{(v, \phi) \in   V \oplus V^* | \, v\otimes \phi \in
[A,\mathrm{gl}(V)] \}.$$ It is convenient to work with $Q_A$ since
its description is linear.

For example, we will show that $Q_{A_1 \oplus A_2} \subset Q_{A_1}
\times Q_{A_2}$. This will allow us to decompose the problem into
Jordan blocks (see section \ref{PfDirectSum}).

Then we will solve the case of one Jordan block (in section
\ref{PfJordanBlock}) using an important result by Rallis and
Schiffmann which is proven using the Weil repesentation.
\section{Preliminaries and notations}
%
%
%
%
%
We will use the standard terminology of $l$-spaces introduced in
\cite{BZ}, section 1. We denote by $\Sc(Z)$ the space of Schwartz
functions on an $l$-space $Z$, and by $\Sc^*(Z)$ the space of
distributions on $Z$ equipped with the weak topology.

We fix a nontrivial additive character $\psi$ of $\F$.
\begin{notation} [Fourier transform]
Let $W$ be a finite dimensional vector space over $\F$. Let  $B$ be
a nondegenerate symmetric bilinear form on $W$. We denote by
$\Fou_B:\Sc^*(W) \to \Sc^*(W)$ the Fourier transform defined using
$B$ and the self-dual measure on $W$.

By abuse of notation, we also denote by $\Fou_B$ the partial Fourier
transform $\Fou_B:\Sc^*(Z \times W) \to \Sc^*(Z\times W)$ for any
$l$-space $Z$.

If $W=U \oplus U^*$ then it has a canonical symmetric bilinear form
given by the quadratic form $Q((v,\phi)):=\langle \phi , v\rangle
:=\phi(v)$. We will denote the Fourier transform defined by it
simply by $\Fou_W$. If there is no ambiguity, we will denote it
simlpy by $\Fou$.
\end{notation}
\begin{proposition} \label{2Four}
Let $W_1 \oplus W_2$ be finite dimensional vector spaces. Let $B_1$
and $B_2$ be nondegenerate symmetric bilinear forms on $W_1$ and
$W_2$ respectively. Let $Z \subset W_1$ be a closed subset. Let $\xi
\in \Sc^*(W_1 \oplus W_2)$ be a distribution. Suppose that
$\Fou_{B_1 \oplus B_2}(\xi)$ is supported in $Z\times W_2$. Then
$\Fou_{B_1}(\xi)$ is also supported in $Z\times W_2$.
\end{proposition}
\emph{Proof. }Let $p_1$ denote the projection $W_1\oplus W_2 \to
W_1$. Since $\Fou_{B_2}$ does not change the projection of the
support of a  distribution to $W_1$,\\
$$p_1(\mathrm{Supp}(\Fou_{B_1}(\xi)))=p_1(\mathrm{Supp}(\Fou_{B_2} \circ
\Fou_{B_1}(\xi)))=p_1(\Fou_{B_1 \oplus B_2}(\xi)) \subset Z$$
\proofend\\
We will use the localization principle, formulated in \cite{Ber},
section 1.4.

\begin{theorem}[Localization principle] \label{LocPrin}
Let  $q:Z \to T$ be a continuous map of $l$-spaces. Denote $Z_t:=
q^{-1}(t)$. Consider $\Sc^*(Z)$ as $\Sc(T)$-module. Let $M$ be a
closed subspace of $\Sc^*(Z)$ which is an $\Sc(T)$-submodule. Then
$M=\overline{\bigoplus_{t \in T} (M \cap \Sc^*(Z_t))}$.
\end{theorem}
Informally, it means that in order to prove a certain property of
distributions on $Z$ it is enough to prove that distributions on
every fiber $Z_t$ have this property.
\begin{corollary} \label{LocPrinCor}
Let $q:Z \to T$ be a continuous map of $l$-spaces. Let an $l$-group
$H$ act on an $l$-space $Z$ preserving the fibers of $q$. Let $\mu$
be a character of $H$. Suppose that for any $t\in T$,
$\Sc^*(q^{-1}(t))^{H,\mu}=0$. Then $\Sc^*(Z)^{H,\mu}=0$
\end{corollary}

\begin{corollary} \label{Product}
Let $H_i \subset \widetilde{H}_i$ be $l$-groups acting on $l$-spaces
$Z_i$ for $i=1, \ldots, k$. Suppose that
$\Sc^*(Z_i)^{H_i}=\Sc^*(Z_i)^{\widetilde{H}_i}$ for all $i$. Then
$\Sc^*(\prod Z_i)^{\prod H_i}=\Sc^*(\prod Z_i)^{\prod
\widetilde{H}_i}$.
\end{corollary}

We will use the following version of the Frobenious reciprocity. It
can be easily deduced from the Frobenious reciprocity described in
\cite{Ber}, section 1.5.

\begin{theorem}[Frobenious reciprocity] \label{Frob}
Let a unimodular $l$-group $H$ act transitively on an $l$-space $Z$.
Let $\varphi:E \to Z$ be an $H$-equivariant map of $l$-spaces. Let
$x\in Z$. Suppose that its stabilizer $\mathrm{Stab}_H(x)$ is
unimodular. Let $W$ be the fiber of $x$.
Let $\mu$ be a character of $H$. Then\\
(i) There exists a canonical isomorphism $\mathrm{Fr}:
\Sc^*(E)^{H,\mu}
\to \Sc^*(W)^{\mathrm{Stab}_H(x),\mu}$.\\
(ii) For any distribution $\xi \in \Sc^*(E)^{H,\mu}$,
$\mathrm{Supp}(\mathrm{Fr}(\xi))=\mathrm{Supp}(\xi)\cap W$.\\
(iii) Frobenious reciprocity commutes with Fourier transform.

Namely, let $W$ be a finite dimensional linear space over $\F$ with
a nondegenerate bilinear form $B$. Let $H$ act on $W$ linearly
preserving $B$.

 Then
for any $\xi \in \Sc^*(Z\times W)^{H,\mu}$, we have
$\Fou_{B}(\mathrm{Fr}(\xi))=\mathrm{Fr}(\Fou_{B}(\xi))$ where
$\mathrm{Fr}$ is taken with respect to the projection $Z \times W
\to Z$.
\end{theorem}
\begin{definition}
Let $W$ be a finite dimensional vector space over $\F$. We call a
distribution $\xi \in \Sc^*(W)$ \textbf{abs-homogeneous of degree
$\mathbf{d}$} if for any function $f \in \Sc(W)$,
$|\xi(h_{t^{-1}}(f))| = |t|^{-d} |\xi(f)|$ where
$(h_{t^{-1}}(f))(v)= f(tv)$.
\end{definition}
 For example, a Haar measure on $W$ is abs-homogeneous of degree
$\dim W$ and the $\delta$-distribution is abs-homogeneous of degree
$0$.\\

A crucial step in the proof of the main theorem is the following
special case of a result by Rallis and Schiffmann (\cite{RS}, lemma
8.1).
\begin{theorem}[Rallis-Schiffmann] \label{Metaplectic}
Let $W$ be a finite dimensional vector space over $\F$ and $B$ be a
nondegenerate symmetric bilinear form on $W$. Denote $Z(B):=\{ v \in
W | B(v,v)=0\}$. Let $\xi$ be a distribution on $W$. Suppose that
both $\xi$ and $\Fou_{B}(\xi)$ are supported in $Z_{B}$.

 Then
$\xi$ is abs-homogeneous of degree $\frac{1}{2} \dim W$.
\end{theorem}
For the benefit of the reader we reproduce the proof of this theorem
in section \ref{PfWeil}.
\begin{remark}
Let $Z$ be an $l$-space and $Q\subset Z$ be a closed subset. We will
identify $\Sc^*(Q)$ with the space of all distributions on $Z$
supported on $Q$. In particular, we can restrict a distribution
$\xi$ to any open subset of the support of $\xi$.
\end{remark}
\section{Reformulations of the problem} \label{Reform}
In this section we will prove the following proposition.

\begin{proposition} \label{reduction}
Theorem \ref{main} implies theorem \ref{Goal}.
\end{proposition}

We will divide this reduction to several propositions.\\
Consider the action of $\tG_n$ on $G_{n+1}$ and on
$\mathrm{gl}_{n+1}$ where $G_n$ acts by conjugation and the
generator of $S_2$ acts by transposition.

\begin{proposition} \label{Red1}
If  $\Sc^*(G_{n+1})^{\tG_n,\chi}=0$ then theorem \ref{Goal} holds.
\end{proposition}
The proof is straightforward.
\begin{proposition} \label{Red2}
If  $\Sc^*(\mathrm{gl}_{n+1})^{\tG_n,\chi}=0$ then
$\Sc^*(G_{n+1})^{\tG_n,\chi}=0$.
\end{proposition}
\emph{Proof.}\footnote{This proof is analogous to the proof of an
analogous statement in \cite{Ber}, section 2.2.} Let $\xi \in
\Sc^*(G_{n+1})^{\tG_n,\chi}$. We have to prove $\xi=0$. Assume the
contrary. Take $p \in \mathrm{Supp}(\xi)$. Let $t=\mathrm{det}(p)$.
Let $f\in \Sc(\F)$ be such that $f(0)=0$ and $f(t) \neq 0$. Consider
the determinant map $\mathrm{det}:G_{n+1} \to \F$. Consider
$\xi':=(f \circ \mathrm{det})\cdot \xi$. It is easy to check that
$\xi' \in \Sc^*(G_{n+1})^{\tG_n,\chi}$ and $p \in
\mathrm{Supp}(\xi')$. However, we can extend $\xi'$ by zero to
$\xi'' \in \Sc^*(\mathrm{gl}_{n+1})^{\tG_n,\chi}$, which is zero by
the assumption. Hence $\xi'$ is also zero. Contradiction. \proofend

\begin{proposition} \label{Red3}
If $\Sc^*(\mathrm{sl}_{n+1})^{\tG_n,\chi}=0$ then
$\Sc^*(\mathrm{gl}_{n+1})^{\tG_n,\chi}=0$
\end{proposition}
\emph{Proof. }Consider the trace map $\mathrm{tr}:\mathrm{gl}_{n+1}
\to \F$. By the localization principle (Corollary \ref{LocPrinCor}),
it is enough to prove that for any $t\in \F$,
$\Sc^*(\mathrm{tr}^{-1}(t))^{\tG_n,\chi}=0$. However, all
$\mathrm{tr}^{-1}(t)$ are isomorphic as $\tG_n$-equivariant
$l$-spaces to $\mathrm{sl}_{n+1}$ by $A \mapsto A -
\frac{\mathrm{tr}(A)}{n}\mathrm{Id}$. \proofend

\begin{proposition} \label{Red4}
If $\Sc^*(X_{n})^{\tG_n,\chi}=0$ then
$\Sc^*(\mathrm{sl}_{n+1})^{\tG_n,\chi}=0$.
\end{proposition}
\emph{Proof. } Consider the map $q:\mathrm{sl}_{n+1} \to \F$ given
by $q(B):= B_{n+1,n+1}$. By the localization principle (corollary
\ref{LocPrinCor}), it is enough to prove that for any $t\in \F$,
$\Sc^*(q^{-1}(t))^{\tG_n,\chi}=0$. However, all $q^{-1}(t)$ are
isomorphic as $\tG_n$-equivariant $l$-spaces to $X_n$ by
$$ \left(
  \begin{array}{cc}
    A_{n \times n} & v_{n\times 1} \\
    \phi_{1\times n} & \lambda \\
  \end{array}
\right) \mapsto (A + \frac{\lambda}{n} \mathrm{Id}, v ,\phi) $$
\proofend

This finishes the proof of proposition \ref{reduction}.
\begin{remark} \label{CorFree}
One can give a coordinate free description of $\tG$ and of its
action on $X$. Namely, $\tG$ is isomorphic to the disjoint union of
the group $G=\mathrm{Aut}(V)$ of automorphisms of $V$ and the set
$\mathrm{Iso}(V,V^*)$ of isomorphisms between $V$ and $V^*$ with the
following multiplication. Let $g,g' \in \mathrm{Aut}(V)$ and $h,h'
\in Iso(V,V^*)$.
$$ g \times g' := g \circ g' \quad\quad h\times g := h \circ g$$
$$ g \times h := {g^*}^{-1} \circ h \quad\quad h\times h' := {h^*}^{-1} \circ h'$$
The action of $\tG$ on $X$ is given by
$$g(A,v,\phi)=(gA,gv,(g^*)^{-1}\phi) \quad \text{ and }
\quad h(A,v,\phi)=((hAh^{-1})^*,(h^*)^{-1} \phi,hv)$$
\end{remark}
\section{Proof of the main theorem}

Rallis and Schiffmann have proven various properties of the support
of $(\tG,\chi)$-equivariant distributions on $X$. We will now
summarize those that we need in our proof.

\begin{notation}
Denote the cone of nilpotent operators in $\mathrm{sl}(V)$ by
$\mathcal{N}$. Denote also
$$Y:=\{(v,\phi)\in V \oplus V^*| \, \langle
\phi,v \rangle=0\}.$$
\end{notation}

\begin{theorem}[Rallis-Schiffmann]\label{ThmRS}
Suppose that the main theorem holds for all dimensions smaller than
$n$ for all finite extensions $\E$ of $\F$. Let $\xi$ be a
$(\tG,\chi)$-equivariant distribution on $X$. Then $\xi$ is
supported in $\mathcal{N} \times Y$.
\end{theorem}
We reproduce the proof of this theorem in the appendix (section
\ref{SecRalSchiff}).

Now we will stratify the nilpotent cone and reduce the support of
the distribution stratum by stratum.

\begin{notation}
For any $i$ we denote by ${\mathcal{N}_i}$ the union of all
nilpotent orbits of dimensions $\leq i$. Note that $\mathcal{N}_i$
are Zariski closed, $\mathcal{N}_{i} = \mathcal{N}$ for $i$ large
enough  and $\mathcal{N}_{-1}=\emptyset$.
\end{notation}

In order to excise the support of the distribution we will use a
family of automorphisms of the problem, which play a crucial role in
our proof.

\begin{notation}
For any  $\lambda \in \F$ we denote by $\nu_{\lambda} :X \to X$ the
homeomorphism defined by
$$\nu_{\lambda}(A,v,\phi):=(A+\lambda v\otimes \phi-\lambda \frac{\langle
\phi,v \rangle}{n}\mathrm{Id},v,\phi).$$
\end{notation}
A simple but important observation is that $\nu_{\lambda}$ commutes
with the action of $\tG$.

\begin{notation}
Let $O$ be a nilpotent orbit of dimension $i$. We set
$${\widetilde{O}} := (O \times Y) \cap \bigcap _{\lambda \in \F} \nu_{\lambda}^{-1}(\mathcal{N}_i \times Y).$$
\end{notation}

To proceed stratum by stratum we will need the following key lemma.

\begin{lemma} [Key lemma]\label{KeyLem}
Let $O$ be a nilpotent orbit. Note that $O\times V \times V^*$ is
$\tG$-invariant. Let $\zeta \in \Sc^*(O\times V \times
V^*)^{\tG,\chi}$. Suppose that $\mathrm{Supp}(\zeta) \subset
\widetilde{O}$ and $\mathrm{Supp}(\Fou_{V \oplus V^*}(\zeta))\subset
\widetilde{O}$. Then $\zeta=0$.
\end{lemma}
The proof will be given in section \ref{SecVVstar} below.

Now we are ready to prove the main theorem.

\begin{theorem}
Any $(\tG,\chi)$-equivariant distribution on $X$ is zero.
\end{theorem}
\emph{Proof}. We prove by downward induction the following
statement: for any $i$, any $\xi \in \Sc^*(X)^{(\tG,\chi)}$ is
supported in $\mathcal{N}_i \times Y$. For $i$ large enough it is
theorem \ref{ThmRS}. Suppose that the statement is true for $i$ and
let us prove it for $i-1$. Let $\xi \in \Sc^*(X)^{\tG,\chi}$. We
need to show that $\xi_{(\mathcal{N}_i \setminus \mathcal{N}_{i-1})
\times Y} = 0$. For this it is enough to show that for any nilpotent
orbit $O$ of dimension $i$, we have $\xi |_{O \times Y}=0$.

Denote $\zeta = \xi |_{O \times V \times V^*}$. We know that
$\mathrm{Supp} (\xi) \subset \mathcal{N}_i \times Y$. On the other
hand, $\nu_{\lambda} (\xi)$ is also $(\tG,\chi)$-equivariant for any
$\lambda$. Therefore $\mathrm{Supp} (\xi) \subset \bigcap \limits
_{\lambda \in \F} \nu_{\lambda}^{-1}(\mathcal{N}_i \times Y)$.

Thus $$\mathrm{Supp}(\zeta) \subset (O \times Y) \cap (\bigcap
_{\lambda \in \F} \nu_{\lambda}^{-1}(\mathcal{N}_i \times Y))=
\widetilde{O}.$$

Since the action of $\tG$ preserves the standard bilinear form on
$V\oplus V^*$, $\Fou_{V \oplus V^*}(\xi)$ is also
$(\tG,\chi)$-equivariant. Note that $\Fou_{V \oplus V^*}(\zeta) =
\Fou_{V \oplus V^*}(\xi)|_{O \times V \times V^*}$ and hence
$\mathrm{Supp}(\Fou_{V \oplus V^*}(\zeta))$ is also contained in
$\widetilde{O}$. Thus by the key lemma $\zeta=0$. \proofend

\section{Proof of the key lemma} \label{SecVVstar}
\begin{notation}
Let $A \in \mathrm{sl}(V)$ be a nilpotent element. Let $O$ be the
orbit of $A$ and $i$ be the dimension of $O$. We denote by ${R_A}$
the fiber at $A$ of the projection $\widetilde{O} \to O$. We
consider $R_A$ as a subset of $V\oplus V^*$.
\end{notation}
Note that $R_A\subset Y$.
\begin{notation}
Let $A \in \mathrm{sl}(V)$ be a nilpotent element. We denote
$${Q_A} :=\{(v, \phi) \in   V \oplus V^* | \,
v\otimes \phi \in [A,\mathrm{gl}(V)] \}.$$
\end{notation}
\begin{lemma}
$R_A \subset Q_A$
\end{lemma}
\emph{Proof. } Let $(v,\phi) \in R_A$. Let $O$ be the orbit of $A$
and $i$ be the dimension of $O$. Consider the Zariski tangent space
$T_A \mathcal{N}_i$ to $\mathcal{N}_i$ at $A$. It coincides with
$T_AO = [A,\mathrm{gl}(V)]$. Since $\langle \phi, v \rangle =0$, we
see that $\mathcal{N}_i$ contains the line $\{A + \lambda v \otimes
\phi\}$.\\ Thus $v\otimes \phi \in T_A
\mathcal{N}_i=[A,\mathrm{gl}(V)]$ and hence $(v,\phi) \in Q_A$.
\proofend
\begin{notation}
Let $A \in \mathrm{sl}(V)$ be a nilpotent element. We denote by
${C_A}$ the stabilizer of $A$ in $G$ and by ${\widetilde{C}_A}$ the
stabilizer of $A$ in $\tG$.
\end{notation}
It is known that $C_A$ is unimodular and hence $\widetilde{C}_A$ is
also unimodular.

The key lemma follows now from Frobenious reciprocity and the
following proposition.
\begin{proposition} \label{SubKey}
Let $A \in \mathrm{sl}(V)$ be a nilpotent element. Let $\eta \in
\Sc^*(V \oplus V^*)^{C_A}$. Suppose that both $\eta$ and
$\Fou(\eta)$ are supported in $Q_A$. Then $\eta \in \Sc^*(V \oplus
V^*)^{\widetilde{C}_A}$.
\end{proposition}
We will call a nilpotent element $A\in \mathrm{sl}(V_k)$ 'nice' if
the previous proposition holds for $A$. Namely, $A$ is 'nice' if any
distribution $\eta \in \Sc^*(V_k \oplus V_k^*)^{C_A}$ such that both
$\eta$ and $\Fou(\eta)$ are supported in $Q_A$ is also
$\widetilde{C}_A$-invariant.

Proposition \ref{SubKey} clearly follows from the following two
lemmas and Jordan decomposition.
\begin{lemma} \label{DirectSum}
Let $A_1 \in \mathrm{sl}(V_k)$ and $A_2 \in \mathrm{sl}(V_l)$ be
nice nilpotent elements. Then $A_1 \oplus A_2 \in
\mathrm{sl}(V_{k+l})$ is nice.
\end{lemma}
\begin{lemma} \label{JordanBlock}
Let $A \in \mathrm{sl}(V_r)$ be a nilpotent Jordan block. Then $A$
is nice.
\end{lemma}
\subsection{Proof of lemma \ref{DirectSum}} \label{PfDirectSum}
We will need the following simple lemma.
\begin{lemma}
$Q_{A_1 \oplus A_2} \subset Q_{A_1} \times Q_{A_2}$.
\end{lemma}
\emph{Proof.} Let $(v,\phi) \in Q_{A_1 \oplus A_2}$. This means that
$v \otimes \phi = [A_1 \oplus A_2,B]$. Let $B = \left(
                                                  \begin{array}{cc}
                                                    B_{11} & B_{12} \\
                                                    B_{21} & B_{22} \\
                                                  \end{array}
                                                \right)$, $v=v_1+v_2$ and $\phi = \phi_1 +\phi_2$  be the
decompositions corresponding to the blocks of $A_1 \oplus A_2$. Then
$v_1 \otimes \phi_1 = [A_1,B_{11}]$ and $v_2 \otimes \phi_2 =
[A_2,B_{22}]$. Hence $(v_1,\phi_1) \in Q_{A_1}$ and $(v_2,\phi_2)
\in Q_{A_2}$. \proofend

Now let us prove lemma \ref{DirectSum}. Let $A_1 \in
\mathrm{sl}(V_{k_1})$ and $A_2 \in \mathrm{sl}(V_{k_2})$ be nice
operators. Let $\eta \in \Sc^*(V_{k_1} \oplus V_{k_2} \oplus
V_{k_1}^* \oplus V_{k_2}^*)^{C_{A_1 \oplus A_2}}$. Suppose that both
$\eta$ and $\Fou_{V_{k_1+k_2}\oplus V_{k_1+k_2}^*}(\eta)$ are
supported in $Q_{A_1} \times Q_{A_2}$. We need to show that $\eta$
is $\widetilde{C}_{A_1 \oplus A_2}$-invariant.
Note that $\widetilde{C}_{A_1}$ acts on $V_{k_1} \oplus V_{k_2}
\oplus V_{k_1}^* \oplus V_{k_2}^*$. Denote
\begin{multline} M:=\{ \alpha
\in \Sc^*(V_{k_1} \oplus V_{k_1}^* \oplus V_{k_2} \oplus
V_{k_2}^*)^{C_{A_1}} |  \text{ both } \alpha \text{ and }
\Fou_{V_{k_1} \oplus V_{k_1}^*}(\alpha)\\ \text{ are supported in }
Q_{A_1} \times V_{k_2} \times V_{k_2}^* \}.\nonumber
\end{multline}
By proposition \ref{2Four}, $\eta \in M$.
The following lemma follows from the fact that $A_1$ is nice using
the localization principle (theorem \ref{LocPrin}).
\begin{lemma}
$M = M^{\widetilde{C}_{A_1}}$.
\end{lemma}

Therefore, $\eta$ is $\widetilde{C}_{A_1}$-invariant. By similar
reasons, $\eta$ is $\widetilde{C}_{A_2}$-invariant. Since $\eta$ is
$C_{A_1 \oplus A_2}$-invariant, we get that $\eta$ is invariant with
respect to $\widetilde{C}_{A_1 \oplus A_2}$.
This completes the proof of lemma \ref{DirectSum}. \proofend
\subsection{Proof of lemma \ref{JordanBlock}} \label{PfJordanBlock}
Let $A \in \mathrm{sl}_{r}$ be the standard nilpotent Jordan block.
\begin{notation}
Denote $${F^i}:=\Ker A^i=\mathrm{Im} A^{r-i} ; \quad
{L^i}:=(F^{r-i})^{\bot} = \mathrm{Im} (A^*)^{r-i} = \Ker (A^*)^i
\subset V_r^*$$
$$\text{ and } {Z}:=\bigcup _{i=0}^r F^i \oplus L^{r-i}.$$
\end{notation}
We will first prove the following lemma from linear algebra.
\begin{lemma} \label{LinAlg}
$Q(A) \subset Z$.
\end{lemma}
\emph{Proof.} Let $(v,\phi) \in Q(A)$. Note that for any $i\geq 0$
and any element $B\in \mathrm{gl}_r, \, \,
\mathrm{tr}(A^i[A,B])=\mathrm{tr}([A,A^iB])=0$. Hence $\langle \phi,
A^iv \rangle =\mathrm{tr}(A^i v \otimes \phi) = 0$. Therefore the
spaces $W: = \Span\{A^iv\}$  and $\Psi := \Span \{(A^*)^i \phi \}$
are orthogonal and thus $\dim W + \dim \Psi \leq r$. Denote $k:=
\dim W$ and $l:=\dim \Psi$.

Consider the set of all non-zero vectors of the form $A^iv$. Since
$A$ is nilpotent, it is easy to see that this set is linearly
independent. Hence $v \in \Ker A^k$ and by the same reasoning $\phi
\in \Ker (A^*)^l$. But since $l \leq r-k$, $\Ker (A^*)^l \subset
\Ker (A^*)^{r-k}$. Hence $(v,\phi) \in F^k \oplus L^{r-k} \subset
Z$. \proofend

Now let $T:V_r \to V_r^*$ be the symmetric nondegenerate bilinear
form defined by $T(e_i)=e^*_{r+1-i}$. By remark \ref{CorFree}, $T$
can be viewed as an element of $\tG_r$. Since $(TAT^{-1})^*=A$,
$T\in \widetilde{C}_A$.

In order to finish the proof of lemma \ref{JordanBlock} it remains
to prove the following lemma.
\begin{lemma} \label{InJordanBlock}
Consider the action of $\F^{\times}$ on $V_r \oplus V_r^*$ defined
by $\rho(\lambda)(v,\phi):=(\lambda v, \lambda^{-1}\phi)$. Let $\eta
\in \Sc^*(V_r \oplus V_r^*)^{\F^{\times}}$. Suppose that
$T(\eta)=-\eta$ and that both $\eta$ and $\Fou(\eta)$ are supported
in $Z$. Then $\eta=0$.
\end{lemma}
\emph{Proof.}\footnote{The ideas of the proof of Lemma
\ref{InJordanBlock} appeared also in \cite{RS}.} We will prove this
lemma by induction on $r$. The case $r=0$ is trivial. Suppose now
that $r \geq 1$ and the lemma is true for all smaller $r$. By
theorem \ref{Metaplectic}, $\eta$ is abs-homogeneous of degree $r$.
Consider $\eta |_{(V_r \oplus V_r^*) \setminus (F^{r-1} \oplus
V_r^*)}$. Since $Z \setminus (F^{r-1} \oplus V_r^*)=(V_r \setminus
F^{r-1}) \oplus \{0\}$, on $Z \setminus (F^{r-1} \oplus V_r^*)$, the
action of $\F^{\times}$ coincides with homothety. Therefore
$\eta|_{(V_r \oplus V_r^*) \setminus (F^{r-1} \oplus V_r^*)}$ is
homothety invariant. On the other hand, it is abs-homogeneous of
degree $r$. Hence it is zero. So $\eta$ is supported in $F^{r-1}
\oplus V_r^*$. By the same reasons $\eta$ is supported in $V_r
\oplus L^{r-1}$. Hence it is supported in $F^{r-1} \oplus L^{r-1}$.

By the same reasoning $\Fou(\eta)$ is supported in $F^{r-1} \oplus
L^{r-1}$. Hence $\eta$ is invariant with respect to translations in
$(F^{r-1} \oplus L^{r-1})^{\bot}$ which is equal to $F^1 \oplus
L^1$. If $r=1$ it implies $\eta=0$. Otherwise it implies that $\eta$
is the pull back of a distribution $\alpha$ on the space $(F^{r-1}
\oplus L^{r-1})/(F^1 \oplus L^1)$ which can be identified with
$V_{r-2} \oplus V_{r-2}^*$.

It is easy to see that $\alpha$ satisfies the conditions of the
lemma for dimension $r-2$. Hence by the induction hypothesis
$\alpha=0$. \proofend
\section{Weil representation and proof of theorem \ref{Metaplectic}}  \label{PfWeil}

The goal of this section is to prove the following theorem.
\begin{theorem}[Rallis-Schiffmann]
Let $W$ be a finite dimensional vector space over $\F$ and $B$ be a
nondegenerate symmetric bilinear form on $W$. Denote $Z(B):=\{ v \in
W | B(v,v)=0\}$. Let $\xi$ be a distribution on $W$. Suppose that
both $\xi$ and $\Fou_{B}(\xi)$ are supported in $Z_{B}$.

 Then
$\xi$ is abs-homogeneous of degree $\frac{1}{2} \dim W$.
\end{theorem}
For the proof we will need the Weil representation.\footnote{We
summarize here those properties of Weil representation that we need.
For a more complete description see for example \cite{Gel}, sections
2.3 and 2.5.}
%
\begin{notation}
Let $t\in \F^{\times}$ be a scalar. We denote
$${a_t}:= \left(
    \begin{array}{cc}
      t & 0 \\
      0 & t^{-1} \\
    \end{array}
  \right), \quad {n_t}:= \left(
    \begin{array}{cc}
      1 & t \\
      0 & 1 \\
    \end{array}
  \right), \quad {\overline{n}_t}:= \left(
    \begin{array}{cc}
      1 & 0 \\
      t & 1 \\
    \end{array}
  \right), \quad {J}:= \left(
    \begin{array}{cc}
      0 & -1 \\
      1 & 0 \\
    \end{array}
  \right). $$
\end{notation}
The following standard lemma follows from Gauss elimination.
\begin{lemma}
The families $n_t$ and $\overline{n}_t$ generate
$\mathrm{SL}(2,\F)$.
\end{lemma}
The following theorem is well known.
\begin{theorem}[Weil representation] \label{Weil}
Let $W$ be a vector space over $\F$ of dimension $d$. Let $B$ be a
symmetric nondegenerate bilinear form on $W$. Then there exists a
projective representation $\pi_B : \mathrm{SL}(2,\F) \to \mathrm{Aut}(\Sc(W))$ such that\\
(i)  for any $g,h \in \mathrm{SL}(2,\F)$, $\pi_B(gh) = u \pi_B(g)
\pi_B(h)$
for some $u \in \C$  such that $|u|=1$\\
(ii)   $\pi_B(\overline{n}_t)(f)(v) = \psi(tB(v,v))f(v)$ \\
(iii)  $\pi_B(a_t)(f)(v) = |t|^{-\frac{d}{2}}f(t^{-1}v)$\\
(iv) $\pi_B(J) = \Fou_B$.
\end{theorem}
We denote the dual representation on $\Sc^*(W)$ by $\pi_{B}^*$.
\subsection{Proof of theorem \ref{Metaplectic}.}
Since $\mathrm{Supp}(\xi) \subset Z(B)$, we have
$\pi_{B}^*(\overline{n}_t)\xi = \xi$ for all $t$. \\ Since
$\mathrm{Supp}(\Fou_{B}(\xi)) \subset Z(B)$,  we have
$$\pi_{B}^*(J^{-1}\overline{n}_tJ)\xi = u_1 \xi \text{ where
} |u_1|=1.$$ Thus $\pi^*_{B}(n_{-t}) \xi = u_1 \xi$. \\
Since the families $n_t$ and $\overline{n}_t$ generate
$\mathrm{SL}(2,\F)$, we have
$$ \pi_{B}^*(a_t)\xi = u_2 \xi \, \text{ where }|u_2|=1.$$ Thus
$\xi$ is abs-homogeneous of degree $\frac{\dim W}{2}$. \proofend
\section{Appendix: Proof of theorem \ref{ThmRS}} \label{SecRalSchiff}
The goal of this section is to prove the following theorem.
\begin{theorem}[Rallis-Schiffmann]
Suppose that the main theorem holds for all dimensions smaller than
$n$ for all finite extensions $\E$ of $\F$. Let $\xi$ be a
$(\tG,\chi)$-equivariant distribution on $X$. Then $\xi$ is
supported in $\mathcal{N} \times Y$.
\end{theorem}

\begin{lemma} \label{Y}
$\mathrm{Supp}(\xi) \subset sl(V) \times Y$.
\end{lemma}

\emph{Proof. }Let $U:=X \setminus  (sl(V)\times Y)$. We have to show
$\Sc^*(U)^{\tG}=0$. Consider the map $p:U \to \F^{\times}$ given by
$p(A,v,\phi)= \langle \phi, v\rangle$. By the localization principle
(corollary \ref{LocPrinCor}), it is enough to show that
$\Sc^*(p^{-1}(\lambda))^{\tG}=0$ for any $\lambda \in \F^{\times}$.

Fix $\lambda \in F^{\times}$. Denote $Z_{\lambda}:=\{(v,\phi) \in V
\oplus V^* | \langle \phi,v \rangle= \lambda\}$. Note that
$Z_{\lambda}$ is a transitive $\tG$-equivariant $l$-space. Define
$pr_2:p^{-1}(\lambda) \to Z_{\lambda}$ by
$pr_2(A,v,\phi):=(v,\phi)$. Note that $pr_2$ is $\tG$-equivariant.
Let $z_0=(e_n,\lambda e_n^*) \in Z_{\lambda}$ where $e_n$ is the
last element of the standard basis of $V_n$ and $e_n^*$ is the last
element of the dual basis of $V_n^*$.

Note that the stabilizer of $z_0$ is isomorphic to $\tG_{n-1}$ and
the fiber $pr_2^{-1}(z_0)$ is isomorphic to $sl(V_n)$ as a
$\tG_{n-1}$-equivariant $l$-space. Hence by Frobenious reciprocity
(theorem \ref{Frob}),
$\Sc^*(p^{-1}(\lambda))^{\tG_n,\chi}=\Sc^*(sl(V_n))^{\tG_{n-1},\chi}$.
By proposition \ref{Red2}, the main theorem for dimension $n-1$
implies that $\Sc^*(sl(V_n))^{\tG_{n-1},\chi}=0$. \proofend

\begin{lemma} \label{semisimple}
Let $A\in sl(V)$ be a non-zero semisimple element. Let $C_A$ be the
stabilizer of $A$ in $G$, and $\widetilde{C}_A$ be the stabilizer of
$A$ in $\tG$. Let $\mathcal{N}_A \subset \mathcal{N}$ be the subset
of all nilpotent operators that commute with $A$. Then
$\Sc^*(\mathcal{N}_A\times V \times V^*)^{C_A}=\Sc^*(\mathcal{N}_A
\times V \times V^*)^{\widetilde{C}_A}$.
\end{lemma}

\emph{Proof. } It is known that a centralizer of a semisimple
element is a Levi subgroup. Hence $C_A$ is isomorphic to $\prod_i
G_{k_i}(\E_i)$ where $k_i < n$ are certain natural numbers and
$\E_i$ are certain finite extensions of $\F$. It is easy to see
(e.g. using remark \ref{CorFree}) that $\widetilde{C}_A$ can be
identified with a subgroup of $\prod_i \tG_{k_i}(\E_i)$.

Therefore, the main theorem for $k_i$ and $\E_i$ and corollary
\ref{Product} of the localization principle imply $\Sc^*(\prod_i
X_{k_i}(\E_i))^{C_A}=\Sc^*(\prod_i
X_{k_i}(\E_i))^{\widetilde{C}_A}$.

The lemma follows now from the fact that $\mathcal{N}_A \times V
\times V^*$ can be identified with a closed subset of $\prod_i
X_{k_i}(\E_i)$ . \proofend

\emph{Proof of theorem \ref{ThmRS}. } By lemma \ref{Y}, $\xi$ is
supported in $sl(V) \times Y$. Hence it is left to show that $\xi$
is supported in $\mathcal{N} \times V \times V^*$.

Let $\mathfrak{X}:=\mathfrak{X}_n$ be the set of all monic
polynomials of degree $n$ in variable $\lambda$. Consider the map
$\Delta:=\Delta_n: X_n \to \mathfrak{X_n}$ that maps $(A,v,\phi)$ to
the characteristic polynomial of $A$. By the localization principle
(corollary \ref{LocPrinCor}) it is enough to show
$\Sc^*(\Delta^{-1}(P))^{\tG,\chi}=0$ for any $P\neq \lambda^n \in
\mathfrak{X_n}$.

Let $\mathcal{R} \subset sl(V)$ be the set of all semisimple
operators with characteristic polynomial $P$. Note that $\tG$ acts
transitively on $\mathcal{R}$. We recall that by the Jordan
decomposition theorem any operator $A$ can by decomposed in a unique
way into a sum of commuting operators $A_s$ and $A_n$ such that
$A_s$ is semi-simple and $A_n$ is nilpotent. This defines a map
$\mathcal{J}:\Delta^{-1}(P) \to \mathcal{R}$ by
$\mathcal{J}(A):=A_s$. It is easy to see that $\mathcal{J}$ is
continuous and $\tG$-equivariant.

The theorem follows now from lemma \ref{semisimple} by Frobenious
reciprocity. \proofend

A. Aizenbud \\
Faculty of Mathematics and Computer Science, Weizmann Institute of
Science POB 26, Rehovot 76100, ISRAEL and Hausdorff Institute for
Mathematics, Bonn. \\
E-mail: aizenr@yahoo.com

$ $\\
D. Gourevitch\\
Faculty of Mathematics and Computer Science, Weizmann Institute of
Science POB 26, Rehovot 76100, ISRAEL and Hausdorff Institute for
Mathematics, Bonn.\\
E-mail: guredim@yahoo.com.
\end{document}